\documentstyle{amsart}
\input{bull-art}
\newcommand{\nc}{\newcommand}
\nc{\rnc}{\renewcommand}

\nc{\f}{\frac}
\nc{\r}{\ref}
\nc{\lt}{\left}
\nc{\rt}{\right}
\nc{\bc}{\begin{center}}
\nc{\ec}{\end{center}}
\nc{\ba}[1]{\begin{array}{#1}}
\nc{\ea}{\end{array}}
\nc{\be}[1]{\begin{equation} \label{#1}}
\nc{\ee}{\end{equation}}
\nc{\bea}{\begin{eqnarray*}}
\nc{\eea}{\end{eqnarray*}}
\nc{\bae}[1]{\begin{eqnarray} \label{#1}}
\nc{\eae}{\end{eqnarray}}
\nc{\nn}{\nonumber \\}
\nc{\al}{\alpha}
\nc{\bt}{\beta}
\nc{\gm}{\gamma}
\nc{\dl}{\delta}
\nc{\ep}{\varepsilon}
\nc{\zt}{\zeta}
\nc{\et}{\eta}
\nc{\th}{\theta}
\nc{\kp}{\kappa}
\nc{\lm}{\lambda}
\nc{\rh}{\rho}
\nc{\sg}{\sigma}
\nc{\ph}{\varphi}
\nc{\ch}{\chi}
\nc{\ps}{\psi}
\nc{\om}{\omega}
\nc{\pr}{^{\prime}}
\nc{\prs}{^{\prime 2}}
\nc{\del}{{\bf \nabla}}
\nc{\lap}{\Delta}
\nc{\dt}{\cdot}
\nc{\ov}{\overline}
\nc{\la}{\langle}
\nc{\ra}{\rangle}
\nc{\gto}{\rightarrow}
\nc{\aeq}{\approx}
\nc{\sbs}{\subset}
\nc{\bi}{\bibitem}

\newcounter{temporary}
\nc{\fn}[1]{\setcounter{temporary}{#1}\fnsymbol{temporary}}
\nc{\rn}[1]{\setcounter{temporary}{#1}\Roman{temporary}}

\nc{\ns}{\normalshape}
\nc{\ints}{\Bbb Z}
\nc{\reals}{\Bbb R}
\nc{\comps}{\Bbb C}
\nc{\rats}{\Bbb Q}
\rnc{\mod}{\bmod\ }

\theoremstyle{definition}
\newtheorem{definition}{Definition}
\newtheorem{example}{Example}
\theoremstyle{remark}
\newtheorem{remark}{Remark} 
\theoremstyle{plain}
\newtheorem{theorem}{Theorem}
\newtheorem{proposition}[theorem]{Proposition}
\newtheorem{lemma}[theorem]{Lemma}

\newtheorem{corollary}[theorem]{Corollary}
\newtheorem{fact}{Fact}

\begin{document}
\def\currentvolume{27}
\def\currentissue{2}
\def\currentyear{1992}
\def\currentmonth{October}
\def\copyrightyear{1992}
\def\currentpages{217-238}

\title[Prevalence: A Translation-Invariant]{PREVALENCE:  A TRANSLATION-INVARIANT 
``ALMOST EVERY''\\
ON INFINITE-DIMENSIONAL SPACES}

\author[B. R. Hunt]{Brian R. Hunt}
\address{Information and Mathematical Sciences Branch 
(Code R44) \\
Naval Surface Warfare Center \\ Silver Spring, Maryland 
20903-5000}
\curraddr{Institute for Physical Science and Technology \\ 
University of
Maryland \\ College Park, Maryland 20742}

\author{Tim Sauer}
\address{Department of Mathematical Sciences \\ George Mason
University \\ Fairfax, Virginia 22030}

\author[J. A. Yorke]{James A. Yorke}
\address{Institute for Physical Science and Technology \\ 
University of
Maryland \\ College Park, Maryland 20742}

\date{Dec. 22, 1991, and in revised form, May 5, 1992}

\subjclass{Primary 28C20, 60B11; Secondary 58F14}

\thanks{The first author was supported by the NSWC 
Independent Research
Program.  The second and third authors were partially 
supported by
the U.S. Department of Energy (Scientific Computing Staff, 
Office of
Energy Research)} 

\maketitle

\begin{abstract}
We present a measure-theoretic condition for a property to 
hold
``almost everywhere'' on an infinite-dimensional vector 
space,
with particular emphasis on function spaces such as $C^k$ 
and $L^p$.
Like the concept of ``Lebesgue almost every'' on 
finite-dimensional
spaces, our notion of ``prevalence'' is translation 
invariant.
Instead of using a specific measure on the entire space, 
we define
prevalence in terms of the class of all probability 
measures with
compact support.  Prevalence is a more appropriate 
condition than the
topological concepts of ``open and dense'' or ``generic'' 
when one
desires a probabilistic result on the likelihood of a 
given property
on a function space.  We give several examples of 
properties which
hold ``almost everywhere'' in the sense of prevalence.  
For instance,
we prove that almost every $C^1$ map on $\reals^n$ has the 
property
that all of its periodic orbits are hyperbolic.
\end{abstract}

\section{Introduction} \label{se1}

Under what conditions should it be said that a given 
property on an
infinite-dimensional vector space is virtually certain to 
hold?  For
example, how are statements of the following type made 
mathematically
precise?

(1) Almost every function $f{:}\  [0,1] \to \reals$ in 
$L^1$ satisfies
$\int_0^1 f(x) dx \neq 0$.

(2)
Almost every sequence $\{a_i\}_{i=1}^\infty$ in $l^2$ has 
the
property that $\sum_{i=1}^\infty a_i$ diverges.

(3)
Almost every $C^1$ function $f{:}\  \reals \to \reals$ has 
the property
that $f'(x) \neq 0$ whenever $f(x) = 0$.

(4)
Almost every continuous function $f{:}\  [0,1] \to \reals$ 
is nowhere
differentiable.

(5)
If $A$ is a compact subset of $\reals^n$ of box-counting 
dimension
$d$, then for $1 \leq k \leq \infty$, almost every $C^k$ 
function $f{:}\ 
\reals^n \to \reals^m$ is one-to-one on $A$, provided that 
$m > 2d$.
(When $A$ is a $C^1$ manifold, the conclusion can be 
strengthened to
say that almost every $f$ is an embedding.)

(6)
If $A$ is a compact subset of $\reals^n$ of Hausdorff 
dimension $d$,
then for $1 \leq k \leq \infty$, almost every $C^k$ 
function $f{:}\ 
\reals^n \to \reals^m$ has the property that the Hausdorff 
dimension
of $f(A)$ is $d$, provided that $m \geq d$.

(7)
For $1 \leq k \leq \infty$, almost every $C^k$ map on 
$\reals^n$ has
the property that all of its fixed points are hyperbolic 
(and further,
that its periodic points of all periods are hyperbolic).

(8) For $4 \leq k \leq \infty$, almost every $C^k$ 
one-parameter family
of vector fields on $\reals^2$ has the property that as 
the parameter is
varied, every Andronov-H\"{o}pf bifurcation which occurs is
``typical'' (in a sense to be made precise later).

In $\reals^n$, there is a generally accepted definition of 
``almost
every'', which is that the set of exceptions has Lebesgue 
measure
zero.  The above statements require a notion of ``almost 
every'' in
infinite-dimensional spaces.  We will be concerned mainly 
with
function spaces such as $L^p$ for $1 \leq p \leq \infty$ 
and $C^k$ for
(integers) $0 \leq k \leq \infty$ on subsets of 
$\reals^n$; many of
these spaces are Banach spaces, and all have a complete 
metric.  The
following are some properties of ``Lebesgue measure zero'' 
and
``Lebesgue almost every'' which we would like to preserve 
on these
spaces.
\begin{enumerate}
\item[(i)]
A measure zero set has no interior (``almost every'' 
implies dense).
\item[(ii)]
Every subset of a measure zero set also has measure zero.
\item[(iii)]
A countable union of measure zero sets also has measure 
zero.
\item[(iv)]
Every translate of a measure zero set also has measure zero.
\end{enumerate}

One could define ``almost every'' on a given function 
space in terms
of a specific measure.  For example, the Wiener measure on 
the
continuous functions is appropriate for some problems.  
However, the
notion of ``almost every'' with respect to such a measure 
violates
property (iv).  The following paragraph illustrates some 
of the
difficulties involved in defining an analogue of Lebesgue 
measure on
function spaces.  We assume all measures to be defined (at 
least) on
the Borel sets of the space.

In an infinite-dimensional, separable\footnote{By {\em 
separable\/} we
mean that the space has a countable dense subset.} Banach 
space, every
translation-invari\-ant measure which is not identically 
zero has the
property that all open sets have infinite measure.  To see 
this,
suppose that for some $\ep$, the open ball of radius $\ep$ 
has finite
measure.  Because the space is infinite dimensional, one 
can construct
an infinite sequence of disjoint open balls of radius 
$\ep/4$ which
are contained in the $\ep$-ball.  Each of these balls has 
the same
measure, and since the sum of their measures is finite, the
$\ep/4$-balls must have measure $0$.  Since the space is 
separable, it
can be covered with a countable collection of 
$\ep/4$-balls, and thus
the whole space must have measure $0$.  (Even if the space 
were not
separable, we would be left with the undesirable property 
that some
open sets have measure zero, violating property (i) above.)

In the absence of a reasonable translation-invariant 
measure on a
given function space, one might hope there is a measure 
which at least
satisfies condition (iv) above; such a measure is called 
{\em
quasi-invariant\/}.  In $\reals^n$, there are an abundance 
of
finite measures which are quasi-invariant, such as 
Gaussian measure.
However, for an infinite-dimensional, locally convex 
topological
vector space, a $\sg$-finite,\footnote{By {\em 
$\sg$-finite\/} we mean
that the entire space can be expressed as a countable 
union of sets of
finite measure.  This rules out measures such as 
``counting measure'',
which assigns to each set its cardinality.} 
quasi-invariant measure
defined on the Borel sets must be identically zero 
\cite{GM,Sud1,Sud2}
(see also pp. 138--143 of \cite{Yam}).

Rather than search for a partial analogue of Lebesgue 
measure on
function spaces, our strategy is to find an alternate 
characterization
of the concepts of ``Lebesgue measure zero'' and 
``Lebesgue almost
every'' which has a natural extension to function spaces.  
Properties
(i)--(iv) alone do not uniquely determine these concepts, 
but there
is a more subtle property which does.  In $\reals^n$, let 
us consider
the class of ``probability measures with compact 
support'', that is,
those measures $\mu$ for which there exists a compact set 
$T \sbs
\reals^n$ such that $\mu(T) = \mu(\reals^n) = 1$.
\begin{enumerate}
\item[(v)]
Let $S$ be a Borel set.  If there exists a probability 
measure $\mu$
with compact support such that every translate of $S$ has
$\mu$-measure zero, then $S$ has Lebesgue measure zero.
\end{enumerate}
Property (v) is proved in \S \r{se2} (see Fact \r{fa6}) by a
simple application of the Tonelli theorem (a variant of 
the Fubini
theorem \cite{DS}).  Notice that conversely, if $S \sbs 
\reals^n$ has
Lebesgue measure zero, then the hypothesis of property (v) 
is
satisfied with $\mu$ equal to (for instance) the uniform 
probability
measure on the unit ball.

Given a probability measure $\mu$ with compact support, we 
can define
a trans\-lation-invariant measure $\tilde{\mu}$ on Borel 
sets $S$ by
$\tilde{\mu}(S) = 0$ if every translate of $S$ has 
$\mu$-measure zero
and $\tilde{\mu}(S) = \infty$ otherwise.  What property 
(v) above says
is that every such measure $\tilde{\mu}$ is greater than 
or equal to
Lebesgue measure on the Borel sets of $\reals^n$.  Thus 
one way to
show that a Borel set is small, in a translation-invariant
probabilistic sense, is to show that $\tilde{\mu}(S) = 0$ 
for some
$\mu$.  Such a strategy is plausible on 
infinite-dimensional spaces
because it is not hard to find probability measures with 
compact
support (for example, uniform measure on a line segment, 
or on the
unit ball of any finite-dimensional subspace).

In general, we will call a Borel set ``shy'' if 
$\tilde{\mu}(S) = 0$
for some probability measure $\mu$ with compact support, 
and we call
any other set shy if it is contained in a shy Borel set 
(just as every
Lebesgue measure zero set is contained in a Lebesgue 
measure zero
Borel set).  We then define a ``prevalent'' set to be a 
set whose
complement is shy.  This definition may not characterize 
all sets for
which the label ``almost every'' is appropriate; our claim 
is rather
that properties which hold on prevalent sets are 
accurately described
as holding ``almost everywhere''.

In  the absence of a probabilistic notion of ``almost 
every'',
statements such as 1--8 above have often been formulated 
in terms of
the topological notion of ``genericity''.  In this 
terminology, a
property on a complete metric space is said to be {\em 
generic\/} if
the set on which it holds is {\em residual\/}, meaning 
that it
contains a countable intersection of open dense 
sets.\footnote{Many
authors require a residual set to be (not just contain) a 
countable
intersection of open dense sets.  Our terminology follows 
\cite{Oxt}.}
The complement of a residual set is said to be of the {\em 
first
category\/}; equivalently, a first category set is a 
countable union
of nowhere dense sets.  The notion of ``first category'' was
introduced by Baire in 1899, and his category theorem 
ensures that a
residual subset of a complete metric space is nonempty, 
and in fact
dense \cite{Oxt}.

The concepts of ``first category'' and ``generic'' have 
formal
similarities to ``measure zero'' and ``almost every'', 
satisfying a
set of properties analogous to (i)--(iv) above.  They also 
agree for
some sets in $\reals^n$; for example, the set of rational 
numbers has
measure zero and is of the first category.  But perhaps 
too much
emphasis has been placed on those examples in which first 
category
sets happen to have measure zero.  Sets which are open and 
dense in
$\reals^n$ can have arbitrarily small Lebesgue measure, 
and residual
sets can have measure zero.

In fact, many properties are known to be topologically 
generic in
$\reals^n$ but have low probability.  While the reader may 
be able to
provide examples from his or her own experience, we 
include some for
completeness.

\begin{example} \ns
For $n \geq 1$ let $U_n = \{x \in [0,1] : 0 < 2^n x 
\pmod{1} <
2^{-n}\}$.  Notice that $V_m = \bigcup_{n > m} U_n$ is 
open and dense
but has measure less than $2^{-m}$.  Hence generically 
points in
$[0,1]$ satisfy $0 < 2^n x \pmod{1} < 2^{-n}$ for 
infinitely many
values of $n$, but the set of such points ($\bigcap_{m 
\geq 1} V_m$) has
measure zero.  A similar construction arises naturally in 
\cite{KY}.
\end{example}

\begin{example} \ns
Here we consider how well real numbers can be approximated 
by
rationals.  The Liouville numbers are the real numbers 
$\lm$ which
have the property that for all $c,n > 0$ there exist 
integers $p$ and
$q > 0$ such that
\[
\lt|\lm - \f{p}{q}\rt| < \f{c}{q^n}.
\]
As in the previous example, the set of Liouville numbers 
is residual
but has Lebesgue measure zero \cite{Oxt}.  In contrast are 
the
Diophantine numbers, real numbers $\mu$ which have the 
property that
for every $\ep >0$ there exists a $c > 0$ such that for 
all integers
$p$ and $q > 0$,
\[
\lt|\mu - \f{p}{q}\rt| > \f{c}{q^{2+\ep}}.
\]
The set of Diophantine numbers is of the first category 
but has full
Lebesgue measure in every interval.
\end{example}

\begin{example} \ns
Arnold studied the family of diffeomorphisms on a circle
\[
f_{\om,\ep}(x) = x + \om + \ep\sin x \pmod{2\pi},
\]
where $0 \leq \om \leq 2\pi$ and $0 \leq \ep < 1$ are 
parameters.
For each $\ep$ we can define the set
\[
S_\ep = \{\om \in [0,2\pi] : f_{\om,\ep} \text{ has a
stable periodic orbit}\}.
\]
For $0 < \ep < 1$, the set $S_\ep$ is a countable union of 
disjoint
open intervals (one for each rational rotation number), 
and is an
open dense subset of $[0, 2\pi]$.  However, the Lebesgue 
measure of
$S_\ep$ approaches zero as $\ep \to 0$.  For small $\ep$, 
the
probability of picking an $\om$ in this open dense set is 
very small.
See pp. 108--109 of \cite{Arn} for more details.
\end{example}

\begin{example} \ns
Consider the dynamics of an analytic map in the complex 
plane near a
neutral fixed point.  Suppose the fixed point is the 
origin; then the
map can be written in the form
\[
z \mapsto e^{2\pi i \al} z + z^2 f(z)
\]
with $0 \leq \al \leq 1$ and $f(z)$ analytic.  Siegel 
\cite{Sie}
proved that for Lebesgue almost every $\al$ (specifically, 
if $\al$
is not a Liouville number), the above map is conjugate to 
a rotation
in a neighborhood of the origin under an analytic change of
coordinates.  On the other hand, Cremer \cite{Cre} 
previously showed
that if $f$ is a polynomial (not identically zero), then 
for a
residual set of $\al$ the above map is not conjugate to a 
rotation in
any neighborhood of the origin.  These results are 
discussed on pp.
98--105 of \cite{Bla}.
\end{example}

\begin{example} \ns
Consider the map $z \mapsto e^z$ on the complex plane.  
Misiurewicz
\cite{Mis} proved that this map is ``topologically 
transitive'', which
implies that a residual set of initial conditions have dense
trajectories.  On the other hand, Lyubich \cite{Lyu1} and 
Rees
\cite{Ree} showed that Lebesgue almost every initial 
condition has a
trajectory whose limit set lies on the real axis (in fact, 
the limit
set is just the trajectory of $0$).  See \cite{Lyu2} for a 
discussion
of both results.
\end{example}

\begin{example} \ns
For many families of dynamical systems in $\reals^2$ 
depending on a
parameter, Newhouse \cite{New} and Robinson \cite{Rob} 
constructed a set
of parameters for which infinitely many attractors 
coexist.  The
constructed set is residual in an interval, but is shown 
in \cite{TY}
and \cite{NT} to have measure zero.
\end{example}

In view of these examples, one might ask why the concept of
``residual'' is used.  Sometimes, one just wants to show 
that a set
obtained by a countable intersection is nonempty, or 
further that it
is dense.  For example, the existence of continuous but 
nowhere
differentiable functions can be proved by showing that 
they form a
residual subset of the continuous functions; this argument 
is due to
Banach (see \S\rn{3}.34.\rn{8} of \cite{Kur}).  Other 
times, one
wants to show that a set is ``large'' in a topological 
sense, perhaps
because there has been no probabilistic alternative.  The 
concept of
``prevalence'' is intended for situations where a 
probabilistic result
is desired.

In \S \r{se2} we formally define prevalence, shyness (the
opposite of prevalence), and related concepts, and develop 
some of
the basic theory of these notions.  Section \r{se3} 
examines the eight
statements from the beginning of this section in the new 
framework.
In \S \r{se4} we develop some of the theory of 
``transversality''
(between functions and manifolds) in the context of 
prevalence, and
use it to prove the third, seventh, and eighth statements. 
 Finally,
\S \r{se5} discusses some other ideas related to prevalence.

\section{Prevalence} \label{se2}

Let $V$ be a complete metric linear space, by which we 
mean a vector
space (real or complex) with a complete metric for which 
addition and
scalar multiplication are continuous.  When we speak of a 
measure on
$V$ we will always mean a nonnegative measure that is 
defined on the
Borel sets of $V$ and is not identically zero.  We write 
$S + v$ for
the translate of a set $S \sbs V$ by a vector $v \in V$.

\begin{definition} \label{df1} \ns
A measure $\mu$ is said to be {\em transverse\/} to a 
Borel set $S
\sbs V$ if the following two conditions hold:
\begin{enumerate}
\item[{\ns (i)}] There exists a compact set $U \sbs V$ for 
which $0 <
\mu(U) < \infty$.
\item[{\ns (ii)}] $\mu(S + v) = 0$ for every $v \in V$.
\end{enumerate}
\end{definition}

Condition (i) ensures that a transverse measure can always 
be
restricted to a finite measure on a compact set (see Fact 
\r{fa2}
below), and in developing the theory of transverse 
measures it is
often useful to think in terms of probability measures 
with compact
support.  For applications it will be convenient to use 
measures which
(like Lebesgue measure) are neither finite nor have 
compact support.
If $V$ is separable, then all measures which take on a 
value other
than $0$ and $\infty$ can be shown to satisfy condition 
(i) \cite{OU}.

\begin{definition} \label{df2} \ns
A Borel set $S \sbs V$ is called {\em shy\/}\footnote{The 
word
``shy'' was suggested to us by J. Milnor.} if there exists 
a measure
transverse to $S$.  More generally, a subset of $V$ is 
called shy if
it is contained in a shy Borel set.  The complement of a 
shy set is
called a {\em prevalent\/} set.
\end{definition}

Strictly speaking, the above concepts could be called 
``translation
shy'' and ``translation prevalent''.  On manifolds for 
which there is
no distinguished set of translations, the corresponding 
theory is more
difficult; this is a topic we do not address in this 
paper.  We again
emphasize that the definitions of shy and prevalent would 
be unchanged
if we required transverse measures to be probability 
measures with
compact support.

Roughly speaking, the less concentrated a measure is, the 
more sets
it is transverse to.  For instance, a point mass is 
transverse to
only the empty set.  Also, we will later show (see Fact 
\r{fa6}) that
if any measure is transverse to a Borel set $S \sbs 
\reals^n$, then
Lebesgue measure is transverse to $S$ too.  When $V$ is 
infinite
dimensional, a convenient choice for a transverse measure 
is often
Lebesgue measure supported on a finite-dimensional
subspace.\footnote{An exact characterization of Lebesgue 
measure on a
given finite-dimensional subspace depends on the choice of 
a basis for
the subspace, but since we are only interested in whether 
or not sets
have measure zero, the choice of basis is unimportant for 
our
purposes.}  For example, Lebesgue measure on the 
one-dimensional space
spanned by a vector $w \in V$ is transverse to a Borel set 
$S \sbs V$
if for all $v \in V$, the set of $\lm \in \reals$ (or 
$\comps$ if $V$
is complex) for which $v + \lm w \in S$ has Lebesgue 
measure zero.  It
immediately follows that every countable set in $V$ is 
shy, and every
proper subspace of $V$ is shy.  Notice that because it is 
possible to
represent an infinite-dimensional space as the continuous 
linear image
of a proper subspace, the continuous linear image of a shy 
set need
not be shy.

We now present some important facts about transversality 
and shyness.
The first follows immediately from the above definitions, 
and in
particular implies that prevalence is translation invariant.
\begin{fact} \label{fa1}
If $S$ is shy, then so is every subset of $S$ and every 
translate of
$S$.
\end{fact}

\begin{fact} \label{fa2}
Every shy Borel set $S$ has a transverse measure which is
finite with compact support.  Furthermore, the support of 
this
measure can be taken to have arbitrarily small diameter.
\end{fact}
\begin{pf}
Let $\mu$ be a measure transverse to a Borel set $S \sbs 
V$.  Then by
condition (i) of Definition \r{df1} it can be restricted 
to a compact
set $U$ of finite and positive measure, and the 
restriction is
certainly also transverse to $S$.  Also, since $U$ is 
compact it can
be covered for each $\ep > 0$ by a finite number of balls 
of radius
$\ep$, and at least one of these balls must intersect $U$ 
in a set of
positive measure.  The intersection of $U$ with the 
closure of this
ball is compact, and the restriction of $\mu$ to this set 
is also
transverse to $S$.
\end{pf}

An immediate consequence of Fact \r{fa2} is that a shy 
Borel set has
no interior.  The same is then true of every shy set, 
since every shy
set is contained in a shy Borel set.  Hence we have the 
following
fact.
\addtocounter{fact}{-1}
\renewcommand{\thefact}{\arabic{fact}$\mbox{}\<^{\prime}$}
\begin{fact} \label{fa2a}
All prevalent sets are dense.
\end{fact}
\renewcommand{\thefact}{\arabic{fact}}

Next, we would like to know that the union of two shy sets 
is also
shy.  Given Borel sets $S,T \sbs V$ containing the 
original sets and
measures $\mu$ and $\nu$ transverse to $S$ and $T$ 
respectively, we
must then find a measure which is transverse to both $S$ 
and $T$.  We
can assume by Fact \r{fa2} that $\mu$ and $\nu$ are finite 
with
compact support.  Then the measure we desire is the 
convolution
$\mu*\nu$ of $\mu$ and $\nu$, defined as follows.
\begin{definition} \label{df3} \ns
Let $\mu$ and $\nu$ be measures on $V$.  Let $\mu \times 
\nu$ be the
product measure of $\mu$ and $\nu$ on the Cartesian 
product $V \times
V$, and for a given Borel set $S \sbs V$ let $S^\Sigma = 
\{(x,y) \in
V \times V : x+y \in S\}$.  Then $S^\Sigma$ is a Borel 
subset of $V
\times V$, and we define $\mu*\nu(S) = \mu \times \nu 
(S^\Sigma)$.
\end{definition}

If $\mu$ and $\nu$ are finite, then $\mu \times \nu$ is 
finite, and
the characteristic function of $S^\Sigma$ is integrable 
with respect
to $\mu \times \nu$.  Then by the Fubini theorem \cite{DS},
\[
\mu*\nu (S) = \int_V \mu(S - y) \,d\nu(y) = \int_V \nu(S - 
x) \,d\mu(x).
\]
We thus have the following fact.
\begin{fact} \label{fa3}
Let $\mu$ and $\nu$ be finite measures with compact 
support.  If
$\mu$ is transverse to a Borel set $S$, then so is
$\mu*\nu$.\footnote{Notice that $\mu*\nu$ has compact 
support because
its support is contained in the continuous image, under 
the mapping
$(x,y) \mapsto x+y$, of the Cartesian product of the 
supports of $\mu$
and $\nu$.}
\end{fact}

From Fact \r{fa3} it follows that the union of two shy 
sets is shy,
and more generally the following fact holds.
\addtocounter{fact}{-1}
\renewcommand{\thefact}{\arabic{fact}$\mbox{}\<^{\prime}$}
\begin{fact} \label{fa3a}
The union of a finite collection of shy sets is shy.
\end{fact}
\renewcommand{\thefact}{\arabic{fact}}

Fact \r{fa3a} extends to countable unions by a slightly more
complicated argument.
\addtocounter{fact}{-1}
\renewcommand{\thefact}{\arabic{fact}$\mbox{}\<^{\prime%
\prime}$}
\begin{fact} \label{fa3b}
The union of a countable collection of shy sets is shy.
\end{fact}
\renewcommand{\thefact}{\arabic{fact}}
\begin{pf}
Given a countable collection of shy subsets of $V$, let 
$S_1, S_2,
\ldots$ be shy Borel sets containing the original sets.  
Let $\mu_1,
\mu_2, \ldots$ be transverse to $S_1, S_2, \ldots$, 
respectively.  By
Fact \r{fa2}, we can assume without loss of generality 
that each
$\mu_n$ is finite and supported on a compact set $U_n$ 
with diameter
at most $2^{-n}$.  By normalizing and translating the 
measures, we
can also assume that $\mu_n(V) = 1$ for all $n$ and that 
each $U_n$
contains the origin.  With 
these assumptions we can define a measure $\mu$ which is 
essentially
the infinite convolution of the $\mu_n$.  We rely on the 
theory of
infinite product measures; see pp. 200--206 of \cite{DS} 
for details.

The infinite Cartesian product $U^\Pi = U_1 \times U_2 
\times \cdots$
is compact by the
Tychonoff theorem \cite{DS} and has a product
measure $\mu^\Pi = \mu_1 \times \mu_2 \times \cdots$ 
defined on its
Borel subsets, with $\mu^\Pi(U^\Pi) = 1$.  Since $V$ is 
complete and
each vector in $U_n$ lies at most $2^{-n}$ away from zero, 
there is a
continuous mapping from $U^\Pi$ into $V$ defined by $(v_1, 
v_2,
\ldots) \mapsto v_1 + v_2 + \cdots$.  The image $U$ of 
$U^\Pi$ under
this mapping is compact, and $\mu^\Pi$ induces a measure 
$\mu$
supported on $U$, given by $\mu(S) = \mu^\Pi(S^\Sigma)$, 
where
$S^\Sigma = \{(v_1,v_2,\ldots) \in U^\Pi : v_1 + v_2 + 
\cdots \in
S\}$.  We will be done if we show that $\mu$ is transverse 
to every
$S_n$.

Since the Cartesian product of measures is associative (and
commutative), we can write $\mu^\Pi = \mu_n \times 
\nu^\Pi_n$ with
$\nu^\Pi_n = \mu_1 \times \cdots \times \mu_{n-1} \times 
\mu_{n+1}
\times \cdots$.  Let $\nu_n$ be the measure on $V$ induced 
by
$\nu^\Pi_n$ under the summation mapping (as $\mu$ was 
induced by
$\mu^\Pi$).  Then $\mu = \mu_n * \nu_n$, and therefore by 
Fact
\r{fa3}, $\mu$ is transverse to $S_n$.  This completes the 
proof.
\end{pf}

We are now in a position to show that the conditions for 
shyness
given in the beginning of this section can be weakened in 
some cases.
First, consider the following definition.
\begin{definition} \label{df4} \ns
A measure $\mu$ is {\em essentially transverse} to a Borel 
set $S
\sbs V$ if condition (i) of Definition \r{df1} holds and 
$\mu(S + v)
= 0$ for a prevalent set of $v \in V$.
\end{definition}

Though essential transversality is weaker than 
transversality, the
following fact holds.
\begin{fact} \label{fa4}
If a Borel set $S \sbs V$ has an essentially transverse 
measure, then
$S$ is shy.
\end{fact}
\begin{pf}
Let $\mu$ be a measure that is essentially transverse to 
$S$.  As in
Fact \r{fa2} we may assume $\mu$ is finite with compact 
support.  The
set of all $v \in V$ for which $\mu(S - v) > 0$ is shy, 
and hence is
contained in a shy Borel set $T$.  Let $\nu$ be a finite 
measure with
compact support which is transverse to $T$.  Then for all 
$w \in V$,
\[
\mu * \nu(S + w) = \int_V \mu(S + w - y) \,d\nu(y) = 0
\]
since the integrand is nonzero only on a subset of $T + w$ 
and $\nu(T
+ w) = 0$.  Thus $\mu * \nu$ is transverse to $S$, and $S$ 
is shy.
\end{pf}

Next let us examine a local definition of shyness and 
prevalence.
\begin{definition} \label{df5} \ns
A set $S \sbs V$ is {\em locally shy\/} if every point in 
the
space $V$ has a neighborhood whose intersection with $S$ is
shy.  A set is {\em locally prevalent\/} if its complement 
is
locally shy.
\end{definition}

Facts \r{fa1}, \r{fa2a}, and \r{fa3a} immediately hold 
also for local
shyness and local prevalence, but whether Fact \r{fa3b} 
does is not
clear in general.  If $V$ is separable though, it turns 
out that the
local definitions of shyness and prevalence are equivalent 
to the
global definitions.  (On the other hand, it is not clear 
that these
notions are the same in spaces such as $L^\infty$ and 
$l^\infty$.)
\begin{fact} \label{fa5}
All shy sets are locally shy.  If $V$ is separable, all 
locally shy
subsets of $V$ are shy.
\end{fact}
\begin{pf}
The first part of this fact is trivial.  To verify the 
second part,
recall that by the Lindel\"{o}f theorem \cite{DS}, if $V$ 
is a
separable metric space then every open cover of $V$ has a 
countable
subcover.  Given a locally shy set $S \sbs V$, the 
neighborhoods whose
intersections with $S$ are shy cover $V$.  Hence by taking 
a countable
subcover, $S$ can be written as a countable union of shy 
sets.  Thus
by Fact \r{fa3b}, $S$ is shy.
\end{pf}

If $V$ is finite dimensional, then shyness and local 
shyness are
equivalent by Fact \r{fa5}.  In this case we can show 
further that
both of these concepts are equivalent to having Lebesgue 
measure
zero.
\begin{fact} \label{fa6}
A set $S \sbs \reals^n$ is shy if and only if it has 
Lebesgue
measure zero.
\end{fact}
\begin{pf}
We need only consider Borel sets, because every Lebesgue 
measure zero
set is contained in a Borel set with measure zero.  If a 
Borel set $S$
has Lebesgue measure zero, then Lebesgue measure is 
transverse to $S$,
and $S$ is shy.  On the other hand, if a Borel set $S$ is 
shy, then
by Fact \r{fa2} there is a finite measure $\mu$ which is 
transverse to
$S$.  Let $\nu$ be Lebesgue measure.  Though $\nu$ is not 
finite, it
is $\sg$-finite, so by the Tonelli theorem \cite{DS} we 
have (as in
the equation preceding Fact \r{fa3}) that
\[
0 = \int_{\reals^n} \mu(S - y) \,d\nu(y) = \int_{\reals^n} 
\nu(S - x)
\,d\mu(x) = \mu(\reals^n) \nu(S).
\]
In other words, $S$ has Lebesgue measure zero.
\end{pf}

Fact \r{fa6} implies that in $\reals^n$, Lebesgue measure 
is a best
possible candidate to be transverse to a given Borel set.  
As we
mentioned earlier, when looking for a transverse measure 
in an
infinite-dimensional space, a useful type of measure to 
try is
Lebesgue measure supported on some finite-dimensional 
subspace.
\begin{definition} \label{df6} \ns
We call a finite-dimensional subspace $P \sbs V$ a {\em 
probe\/} for
a set $T \sbs V$ if Lebesgue measure supported on $P$ is 
transverse
to a Borel set which contains the complement of $T$.
\end{definition}

Then a sufficient (but not necessary) condition for $T$ to 
be
prevalent is for it to have a probe.  One advantage of 
using probes
is that a single probe can often be used to show that a 
given property
is prevalent on many different function spaces by applying 
the
following simple fact.
\begin{fact} \label{fa7}
If $\mu$ is transverse to $S \sbs V$ and the support of 
$\mu$ is
contained in a subspace $W \sbs V$, then $S \cap W$ is a shy
subset of $W$.
\end{fact}

Next we use one-dimensional probes to show that all 
compact subsets of
an infinite-dimensional space are shy.  We prove in fact 
that given a
compact set $S \sbs V$, there are one-dimensional 
subspaces $L$ for
which every translate of $L$ intersects $S$ in at most one 
point.  To
do this we show that a residual set of vectors in $V$ span
one-dimensional subspaces $L$ with the above property.  
Here then is
an application of the fact that a residual set is nonempty.
\begin{fact} \label{fa8}
If $V$ is infinite dimensional, all compact subsets of $V$ 
are shy.
\end{fact}
\begin{pf}
We assume $V$ is a real vector space; the proof is nearly 
identical
for a complex vector space.  Let $S \sbs V$ be a compact 
set, and
define the function $f : \reals \times S \times S \to V$ by
\[
f(\al, x, y) = \al(x - y).
\]
If a vector $v \in V$ is not in the range of $f$, then $v$ 
spans a
line $L$ such that every translate of $L$ meets $S$ in at 
most one
point; in particular, $L$ is a probe for the complement of 
$S$.  We
then need only show that the range of $f$ is not all of 
$V$; we show
in fact that it is a first category set.  For each 
positive integer
$N$, the set $[-N,N] \times S \times S$ is compact, and 
hence so is
its image under $f$.  Thus the range of $f$ is a countable 
union of
compact sets.  Since $V$ is infinite dimensional, a 
compact set in $V$
has no interior (see p.\ 23 of \cite{Sch}), and is then 
nowhere dense
(because it is closed).  Therefore the range of $f$ is of 
first
category as claimed.
\end{pf}

\section{Applications of prevalence} \label{se3}

From now on, when we say ``almost every'' element of $V$ 
satisfies a
given property, we mean that the subset of $V$ on which 
the property
holds is prevalent.  Given this terminology, the eight 
numbered
statements from the introduction can be proved by 
constructing
appropriate probes (see Definition \r{df6}).

\begin{proposition} \label{pr1}
Almost every function $f: [0,1] \to \reals$ in $L^1$ 
satisfies
$\int_0^1 f(x)\, dx 
\neq 0$.
\end{proposition}

A probe for Proposition \r{pr1} is the one-dimensional 
space of all
constant functions.  Notice that this probe is contained 
in $C^k$ for
$0 \leq k \leq \infty$, so the above property also holds 
for almost
every $f$ in $C^k$.  Similar remarks can be made about 
most of the
results below.

\begin{proposition} \label{pr2}
For $1 < p \leq \infty$, almost every sequence 
$\{a_i\}_{i=1}^\infty$
in $l^p$ has the property that $\sum_{i=1}^\infty a_i$ 
diverges.
\end{proposition}

For Proposition \r{pr2}, the one-dimensional space spanned 
by the
element $\{1/i\}_{i=1}^\infty
 \in l^p$ is a probe for each $1 < p
\leq \infty$.

The third statement in the introduction can be proved using 
the space of constant functions as a probe; this follows 
immediately
from the Sard theorem \cite{Sar}.  Here we state a more 
general
result, which uses a higher-dimensional probe.  We write 
$f^{(i)}$ for
the $i$th derivative of $f$.
\begin{proposition} \label{pr3}
Let $k$ be a positive integer.  Almost every $C^k$ 
function $f{:}
\ \reals \to \reals$ has the property that at each $x \in 
\reals$, at
most one of $\{f^{(i)}(x) : 0 \leq i \leq k \}$ is zero.
\end{proposition}

The space of polynomials of degree $\leq k$ is a probe for 
Proposition
\r{pr3}, as we will prove in the next section.  By Fact 
\r{fa3b},
Proposition \r{pr3} has the following corollary.
\addtocounter{theorem}{-1}
\renewcommand{\thetheorem}{\arabic{theorem}$\mbox{}\<^{%
\prime}$}
\begin{proposition} \label{pr3a}
Almost every $C^\infty$ function $f{:}\  \reals \to 
\reals$ has the
property that at each $x \in \reals$, at most one of 
$\{f^{(i)}(x) :
i \geq 0\}$ is zero.
\end{proposition}
\renewcommand{\thetheorem}{\arabic{theorem}}

Because the dimension of the probe used to prove Proposition
\r{pr3} goes to infinity as $k \gto \infty$, it is not clear
whether Proposition \r{pr3a} can be proved directly using 
a probe.

\begin{proposition} \label{pr4}
Almost every continuous function $f{:}\ [0,1] \to \reals$ 
is nowhere
differentiable.
\end{proposition}

Proposition \r{pr4} requires a two-dimensional probe.  A
one-dimensional probe would be spanned by a continuous 
function $g$
with the property that for all continuous $f{:}\  [0,1] 
\to \reals$, the
function $f + \lm g$ is nowhere differentiable for almost 
every $\lm
\in \reals$.  But if $f(x) = -xg(x)$, then $f + \lm g$ is
differentiable at $x = \lm$ for every $\lm$ between 0 and 
1.  The
proof of Proposition \r{pr4} uses a probe spanned by a 
pair of
functions $g$ and $h$ for which $\lm g + \mu h$ is nowhere
differentiable for every $(\lm,\mu) \in \reals^2$ aside 
from the
origin \cite{Hun}.

Next we state a prevalence version of the Whitney 
Embedding Theorem.
\begin{proposition} \label{pr5}
Let $A$ be a compact $C^1$ manifold of dimension $d$ 
contained in
$\reals^n$.  For $1 \leq k \leq \infty$, almost every 
$C^k$ function
$f {:}\  \reals^n \to \reals^{2d+1}$ is an embedding of $A$.
\end{proposition}

The probe used in the proof of Proposition \r{pr5} is the 
space of
linear functions from $\reals^n$ to $\reals^{2d+1}$.  
Whitney \cite{Whi} showed
that a residual subset of the $C^k$ functions from 
$\reals^n$ to
$\reals^{2d+1}$ are embeddings of $A$.  This result was 
preceded by
a topological version due to Menger in 1926 (see p. 56 of 
\cite{HW}),
which states that for a compact space $A$ of topological 
dimension
$d$, a residual subset of the continuous functions from 
$A$ to
$\reals^{2d+1}$ are one-to-one.  Proposition \r{pr5}, and 
the following
generalization to compact subsets of $\reals^n$ which need 
not be manifolds
(or even have integer dimension), are proved in \cite{SYC}.
\addtocounter{theorem}{-1}
\renewcommand{\thetheorem}{\arabic{theorem}$\mbox{}\<^{%
\prime}$}
\begin{proposition} \label{pr5a}
If $A$ is a compact subset of $\reals^n$ of box-counting 
\RM(capacity\RM)
dimension $d$, and $1 \leq k \leq \infty$, then almost 
every $C^k$
function $f{:}\  \reals^n \to \reals^m$ is one-to-one on 
$A$, provided
that $m > 2d$.
\end{proposition}
\renewcommand{\thetheorem}{\arabic{theorem}}

Our next proposition concerns the preservation of Hausdorff
dimension under smooth transformations.  Once again the 
probe is the
space of all linear functions from $\reals^n$ to 
$\reals^m$; see
\cite{SY} for a proof.
\begin{proposition} \label{pr6}
If $A$ is a compact subset of $\reals^n$ of Hausdorff 
dimension $d$,
and $1 \leq k \leq \infty$, then for almost every $C^k$ 
function $f{:}\ 
\reals^n \to \reals^m$ the Hausdorff dimension of $f(A)$ 
is $d$,
provided that $m \geq d$.
\end{proposition}

\begin{remark} \ns
It is interesting that Proposition \ref{pr5a} fails for 
Hausdorff
dimension (see \cite{SYC}), and Proposition \ref{pr6} 
fails for
box-counting dimension (see \cite{SY}), under any reasonable
definition of ``almost every''.
\end{remark}

We now present a result about the prevalence of 
hyperbolicity for
periodic orbits of maps.  We say that a period $p$ point 
of a map $f{:}\ 
\reals^n \to \reals^n$ is {\em hyperbolic\/} if the 
derivative of the
$p$th iterate of $F$ at the point has no eigenvalues (real 
or complex)
with absolute value $1$.
\begin{proposition} \label{pr7}
Let $p$ be a positive integer.  For $1 \leq k \leq 
\infty$, almost
every $C^k$ map on $\reals^n$ has the property that all of 
its
periodic points of period $p$ are hyperbolic.
\end{proposition}

Proposition \r{pr7} is proved in the next section using 
the space of
polynomial functions of degree at most $2p-1$ as a probe.
Proposition \r{pr7} and Fact \r{fa3b} imply the following 
more
elegant result.
\addtocounter{theorem}{-1}
\renewcommand{\thetheorem}{\arabic{theorem}$\mbox{}\<^{%
\prime}$}
\begin{proposition} \label{pr7a}
For $1 \leq k \leq \infty$, almost every $C^k$ map on 
$\reals^n$ has
the property that all of its periodic points are hyperbolic.
\end{proposition}
\renewcommand{\thetheorem}{\arabic{theorem}}

Next consider one-parameter families of dynamical systems. 
 As the
parameter varies, it is likely that nonhyperbolic points 
will be
encountered, and at such points bifurcations (creation or 
destruction
of periodic orbits, or changes in stability of orbits) can 
occur.  In
general one can expect to prove results of the sort that 
for dynamical
systems of a given type, almost every one-parameter family 
has the
property that all of its bifurcations are 
``nondegenerate'' in some
fashion.  A complete discussion of such results is beyond 
the scope of
this paper, but we include as an example a result about
Andronov-H\"{o}pf bifurcations for flows in the plane.  
For flows (as
opposed to maps), a fixed point is hyperbolic if the 
linear part of
the vector field at the fixed point has no eigenvalues on 
the
imaginary axis.  Generally, a zero eigenvalue results in a 
saddle-node
bifurcation and a pair of nonzero, pure imaginary 
eigenvalues results
in an Andronov-H\"{o}pf bifurcation; see \cite{GH} for 
details.  The
following proposition is proved in the next section.
\begin{proposition} \label{pr8}
For $4 \leq k \leq \infty$, almost every $C^k$ 
one-parameter family
of vector fields $f(\mu,x) {:}\  \reals \times \reals^2 
\to \reals^2$ has
the property that whenever $f(\mu_0,x_0) = 0$ and $D_x 
f(\mu_0,x_0)$
has nonzero, pure imaginary eigenvalues, the flow $\dot{x} 
= f(\mu,x)$
undergoes a nondegenerate Andronov-H\"{o}pf bifurcation in 
the sense
that the following conditions hold in a neighborhood $U$ of
$(\mu_0,x_0)$\RM:
\begin{enumerate}
\item[{\ns (i)}] The fixed points in $U$ form a curve 
$(\mu,x(\mu))$,
  where $x(\mu)$ is $C^k$.
\item[{\ns (ii)}] The point $(\mu,x(\mu))$ is attracting 
when $\mu$
  is on one side of $\mu_0$ and repelling when $\mu$ is on 
the other
  side.
\item[{\ns (iii)}] There is a $C^{k-2}$ 
surface\footnote{The surface
  is proved to be $C^{k-2}$ in \cite{MM}.  However, we 
suspect that
  this surface can actually be shown to be $C^{k-1}$, in 
which case
  this proposition applies to $C^3$ vector fields also.} 
of periodic
  orbits in $\reals \times \reals^2$ which has a quadratic 
tangency
  with the plane $\mu = \mu_0$.  The periodic orbits are 
attracting if
  the fixed points for the same parameter values are 
repelling, and
  are repelling if the corresponding fixed points are 
attracting.
\end{enumerate}
\end{proposition}

\section{Transversality and prevalence} \label{se4}

The proofs of Propositions \r{pr3}, \r{pr7}, and \r{pr8} 
are based on
the idea of ``transversality'', which we will discuss now 
in the
context of functions from one Euclidean space to another.  
Given $1
\leq k \leq \infty$ and $0 \leq d < \infty$, we call $M 
\sbs \reals^n$
a {\em $C^k$ manifold\/} of dimension $d$ if for all $x 
\in M$ there
is an open neighborhood $U \sbs \reals^n$ of $x$ and a $C^k$
diffeomorphism $\ph {:}\  U \to V \sbs \reals^n$ such that 
$\ph(M \cap U)
= (\reals^d \times \{0\}) \cap V$.  The {\em tangent 
space\/} to $M$
at $x$, denoted by $T_x M$, is defined to be the inverse 
image of
$\reals^d \times \{0\}$ under $D\ph(x)$.  Notice that an 
open set in
$\reals^n$ is a $C^\infty$ manifold of dimension $n$, with 
tangent
space $\reals^n$ at every point.
\begin{definition} \label{df7} \ns
Let $A \sbs \reals^n$ and $Z \sbs \reals^m$ be manifolds.  
We say that
a $C^1$ function $F {:}\  A \to \reals^m$ is {\em 
transversal\/} to $Z$ if
whenever $F(x) \in Z$, the spaces $DF(x)(T_x A)$ and 
$T_{F(x)} Z$ span
$\reals^m$.
\end{definition}
\begin{remark} \ns
If $DF(x)$ maps $T_x A$ onto $\reals^m$ for all $x \in A$, 
then $F$ is
transversal to every manifold in $\reals^m$; in this case 
we say that
$F$ is a {\em submersion\/}.
\end{remark}

In our applications $A$ is always an open set in 
$\reals^n$, so
the results below are stated only for this case, though 
they remain
valid for functions whose domains are sufficiently smooth 
manifolds.
A basic result is the following (see \cite{GP} for a proof).
\addtocounter{theorem}{-8}
\begin{theorem}[Parametric Transversality Theorem] 
\label{th1}
Let $B \sbs \reals^q$ and $A \sbs \reals^n$ be open sets.  
Let $F {:}\  B
\times A \to \reals^m$ be $C^k$, and let $Z$ be a $C^k$ 
manifold of
dimension $d$ in $\reals^m$, where $k > \max(n+d-m,0)$.  
If $F$ is
transversal to $Z$, then for almost every $\lm \in B$, the 
function
$F(\lm,\cdot) : A \to \reals^m$ is transversal to $Z$.
\end{theorem}

Notice that if $F {:}\  A \to \reals^m$ is transversal to 
$Z \sbs
\reals^m$ and the codimension of $Z$ (that is, $m$ minus 
the dimension
of $Z$) is greater than the dimension of $A$, then $F(A)$ 
cannot
intersect $Z$.  This observation is the basis for the 
following
general scheme for proving results like Propositions 
\r{pr3}, \r{pr7},
and \r{pr8}.  To show that almost every $f$ in a space 
such as
$C^k(\reals^n)$ has a given property, construct a function 
$F$
consisting of the derivatives of $f$ up to a certain 
order, and let
$Z$ be a manifold defined by a set of $n+1$ conditions 
which $F$ must
satisfy at some point in $\reals^n$ in order for $f$ not 
to have the
desired property.  By an appropriate generalization of 
Theorem
\r{th1}, it will follow that for almost every $f$, there 
is no point
in $\reals^n$ at which $F$ satisfies the undesirable 
conditions.

Let us formalize the above procedure.
\begin{definition} \label{df8} \ns
Let $A \sbs \reals^n$ be open, and let $f {:}\  A \to 
\reals^m$ be
$C^l$.  For $k \leq l$, we define the {\em $k$-jet\/} of 
$f$ at
$x$, denoted $j^k f(x)$, to be the ordered pair consisting 
of $x$ and
the degree $k$ Taylor polynomial of $f$ at $x$.  Then $j^k 
f$ is a
$C^{l - k}$ function from $A$ to a space $J^k(\reals^n, 
\reals^m) =
\reals^n \times P^k(\reals^n, \reals^m)$, where 
$P^k(\reals^n,
\reals^m)$ is the space of polynomials of degree $\leq k$ 
from
$\reals^n$ to $\reals^m$.  We write
\[j^k f(x) = (x, f(x), Df(x), \ldots, D^k f(x)),\]
where the coordinates $(f(x), Df(x), \ldots, D^k f(x))$ 
represent the
(unique)
polynomial in
$P^k(\reals^n, \reals^m)$ with the same derivatives up
to order $k$ as $f$ at $x$.
\end{definition}
\begin{remark} \ns
We will later write $J^k(\reals^n, \reals^m) = 
J^{k-1}(\reals^n,
\reals^m) \times \widehat{P}^k(\reals^n, \reals^m)$, where
$\widehat{P}^k(\reals^n, \reals^m)$ can be thought of as 
the space of
homogeneous degree $k$ polynomials from $\reals^n$ to 
$\reals^m$.
More precisely, $j^k f(x)$ can be decomposed into 
$(j^{k-1} f(x), D^k
f(x))$, where $D^k f(x)$ represents a degree $k$ 
polynomial which is
homogeneous in a coordinate system based at $x$.

The following is an example of the type of result we need; 
it is a
prevalence version of a result previously formulated in 
terms of
genericity \cite{Hir}.
\end{remark}
\addtocounter{theorem}{-1}
\renewcommand{\thetheorem}{\arabic{theorem}$\mbox{}\<^{%
\prime}$}
\begin{theorem}[Jet Transversality Theorem] \label{th1a}
Let $A \sbs \reals^n$ be open and let $Z$ be a $C^r$ 
manifold in
$J^k(\reals^n, \reals^m)$ with codimension $c$, where $r >
\max(n-c,0)$.  For $k + \max(n-c,0) < l \leq \infty$, 
almost every
$C^l$ function $f : A \to \reals^m$ has the property that 
$j^k f$
is transversal to $Z$.
\end{theorem}
\renewcommand{\thetheorem}{\arabic{theorem}}
\begin{pf}
Let $P = P^k(\reals^n, \reals^m)$, thinking of $P$ for now 
as a
subspace of the $C^l$ functions from $A$ to $\reals^m$.  
We claim
that $P$ is a probe (see Definition \r{df6}) for the above 
property.
For $p \in P$, let $f_p(x) = f(x) + p(x)$, and define the 
function $F
{:}\  P \times A \to J^k(\reals^n, \reals^m)$ by $F(p,x) = 
j^k f_p(x)$.
Notice that $F$ is a submersion, because the first $n$ 
coordinates of
$F$ are just $x$, and for a given $x$ the remaining 
coordinates of $F$
act as a translation (by the Taylor polynomial of $f$ at 
$x$) on $P$.
In particular, $F$ is transversal to $Z$.  Then by Theorem 
\r{th1},
for almost every $p \in P$, the function $F(p,\cdot) = j^k 
f_p$ is
transversal to $Z$, and therefore $P$ is a probe as claimed.
\end{pf}

A special case of Theorem \r{th1a} is the following 
prevalence
version of the Thom Transversality Theorem.
\addtocounter{theorem}{-1}
\renewcommand{\thetheorem}{\arabic{theorem}$\mbox{}\<^{%
\prime\prime}$}
\begin{corollary} \label{co1b}
Let $A \sbs \reals^n$ be open and let $Z$ be a $C^r$ 
manifold in
$\reals^m$ with codimension $c$, where $r > \max(n-c,0)$.  
For
$\max(n-c,0) < k \leq \infty$, almost every $C^k$ function 
$f {:}\  A \to
\reals^m$ is transversal to $Z$.
\end{corollary}
\renewcommand{\thetheorem}{\arabic{theorem}}

Propositions \r{pr3}, \r{pr7}, and \r{pr8} can be proved 
using
Theorems \r{th1} and \r{th1a}, except that we would then 
have to assume
in Proposition \r{pr3} that $f$ is $C^{k+1}$ and in 
Proposition
\r{pr7} that the map is $C^2$.  Instead we use the 
following results,
which do
not require those additional assumptions and also allow us
to avoid determining the entire manifold structure of $Z$.
\begin{definition} \label{df9} \ns
We say that a set $S$ is a {\em zero set\/} in a manifold 
$M$ of
dimension $d$ if $S \sbs M$ and for every $x \in M$ there 
is a
neighborhood $U$ of $x$ and a diffeomorphism $\ph$ on $U$ 
which takes
$M \cap U$ to an open set in $\reals^d$ and for which 
$\ph(S \cap U)$
has Lebesgue measure zero in $\reals^d$.
\end{definition}
\begin{remark} \ns
Since sets of Lebesgue measure zero are preserved under
diffeomorphism, the particular choice of $\ph$ in 
Definition \r{df9}
does not matter; that is, a zero set in $M$ has Lebesgue 
measure zero
with respect to all local $C^1$ coordinate systems on $M$.

We will need a Fubini-like result for zero sets of 
manifolds which
allows us to prove that a Borel set is a zero set in $M$ 
by showing
that it is a zero set on the leaves of an appropriate 
foliation of
$M$.  See \cite{PS} for a general result of this type; for 
our
purposes we need only the following simple lemma, which 
follows
directly from the Fubini theorem.
\end{remark}
\begin{lemma} \label{le2}
Let $M$ be a manifold of dimension $d$, and let 
$\{M_\al\}$ be a
partition of $M$ into manifolds of dimension $d\pr < d$ 
with the
following property\,\RM:  every $x \in M$ has a 
neighborhood $U$ and a
diffeomorphism $\ph$ on $U$ which maps $M \cap U$ to an 
open set in
$\reals^d$ and which maps those $M_\al$ which intersect 
$U$ into
parallel hyperplanes of dimension $d\pr$.  If $Z$ is a 
Borel set
in $M$ and $Z \cap M_\al$ is a zero set in $M_\al$ for 
every $\al$,
then $Z$ is a zero set in $M$.
\end{lemma}
\begin{remark} \ns
The hypotheses of Lemma \r{le2} are satisfied if $M$ can 
be written as
$M_1 \times M_2$ with $M_1$, $M_2$ manifolds, and the 
partition of $M$
consists of all manifolds of the form $\{x\} \times M_2$ 
with $x \in
M_1$; this will usually be the case when we apply Lemma 
\r{le2}.
\end{remark}

We now present measure-theoretic analogues to Theorems 
\r{th1} and
\r{th1a}.
\begin{lemma}[Measure Transversality Lemma] \label{le3}
Let $B \sbs \reals^q$ and $A \sbs \reals^n$ be open sets, 
with points
in $B$ denoted by $\lm$ and points in $A$ denoted by $x$.  
Let $F {:}\  B
\times A \to \reals^m \times \reals^s$ be a continuous 
function
with components $G {:}\  B \times A \to \reals^m$ and $H 
{:}\  B \times A \to
\reals^s$.  Assume that the derivatives $D_\lm G$, $D_x 
G$, and
$D_\lm H$ exist and are continuous at every point of $B 
\times A$ \RM(but
$D_x H$ need not exist\RM).  Let $M$ be a manifold in 
$\reals^m$ with
codimension $n$, and assume that for all $x \in A$ and all 
$y \in
\reals^s$, the function $F(\cdot,x)$ is transversal to $M 
\times
\{y\}$.  Let $Z$ be a zero set in $M \times \reals^s$.  
Then for
almost every $\lm \in B$, there is no $x \in A$ for which 
$F(\lm,x)
\in Z$.
\end{lemma}
\begin{remark} \ns
The transversality hypothesis of Lemma \r{le3} is 
automatically
satisfied if $D_\lm F$ has full rank at every point of $B 
\times A$
(that is, if $F(\cdot, x)$ is a submersion for every $x 
\in A$).
\end{remark}

Lemma \r{le3} will be proved at the end of this section.  
Notice that
in the case that $F$ is $C^1$ and $Z$ is a manifold with 
codimension
greater than $n$, Lemma \r{le3} is a special case of 
Theorem \r{th1}.
In much the same way as Theorem \r{th1a} followed from 
Theorem
\r{th1}, the next theorem follows from Lemma \r{le3}.
\addtocounter{theorem}{-1}
\renewcommand{\thetheorem}{\arabic{theorem}$\mbox{}\<^{%
\prime}$}
\begin{theorem}[Measure Jet Transversality Theorem] 
\label{th3a}
Assume $k \geq 1$.  Let $A \sbs \reals^n$ be open and let 
$M$ be a
manifold in $J^{k-1}(\reals^n, \reals^m)$ with codimension 
$n$.  Let
$\pi$ be the projection from $J^{k-1}(\reals^n, \reals^m)$ 
onto its
first $n$ coordinates, and assume that $\pi|_M$ is a 
submersion.  Let
$Z$ be a zero set in $M \times \widehat{P}^k(\reals^n, 
\reals^m)$ \RM(see the
remark following Definition \r{df8}\RM).  Then for $k \leq 
l \leq
\infty$, almost every $C^l$ function $f : A \to \reals^m$ 
has the
property that the image of $A$ under $j^k f$ does not 
intersect $Z$.
\end{theorem}
\renewcommand{\thetheorem}{\arabic{theorem}}
\begin{remark} \ns
In our applications, $M$ will be defined by a set of 
conditions that
$f$ and its derivatives must satisfy at some point $x$.  
When these
conditions do not explicitly depend on $x$, the hypothesis 
that
$\pi|_M$ be a submersion is trivially satisfied.
\end{remark}
\begin{pf}
The proof is the same as for Theorem \r{th1a}, except that 
we must
verify that $F(p,x) = j^k (f(x) + p(x))$, where $p \in 
P^k(\reals^n,
\reals^m)$, satisfies the transversality condition of 
Lemma \r{le3}.
Given $y \in \widehat{P}^k(\reals^n, \reals^m)$, we have 
by hypothesis
that under projection onto the first $n$ coordinates in 
$J^k(\reals^n,
\reals^m)$, the tangent space to $M \times \{y\}$ at any 
point
projects onto all of $\reals^n$.  The remaining 
coordinates in
$J^k(\reals^n, \reals^m)$ are just $P^k(\reals^n, 
\reals^m)$,
and when composed with projection onto the latter space, 
$F(\cdot,x)$
is just a translation (and hence a submersion) for every 
$x$.  Thus
$F(\cdot,x)$ is transversal to $M \times \{y\}$, and the 
hypotheses of
Lemma \r{le3} are satisfied.
\end{pf}

Theorem \r{th3a} says, roughly speaking, that given $n$
``codimension one'' conditions on the $(k-1)$-jet of $f$ 
(these
conditions can depend on $x$, but none can depend only on 
$x$) and an
additional ``measure zero'' condition on the $k$-jet of 
$f$, almost
every $C^k$ function $f$ on a given set in $\reals^n$ does 
not satisfy
all $n+1$ conditions at any point in its domain.  We now 
use Lemma
\r{le3} and Theorem \r{th3a} to prove Propositions 
\r{pr3}, \r{pr7},
and \r{pr8}.

\begin{pf*}{Proof of Proposition \r{pr3}}
We will show for each pair $(i_1,i_2)$ with $0 \leq i_1 < 
i_2 \leq k$
that almost every $f$ in $C^k(\reals)$ has the property that
$f^{(i_1)}(x)$ and $f^{(i_2)}(x)$ are never both zero.  
Let $M$ be the
manifold in $J^{k-1}(\reals,\reals)$ defined by $f^{(i_1)} 
= 0$.  Then
$M$ has codimension $1$, and the set $Z \sbs J^k(\reals, 
\reals)$
defined by $f^{(i_1)} = f^{(i_2)} = 0$ is a zero set in $M 
\times
\widehat{P}^k(\reals, \reals)$.  Therefore by Theorem 
\r{th3a}, almost
every $f$ in $C^k(\reals)$ has the property that $j^k 
f(x)$ is not in
$Z$ for any $x \in \reals$, which is exactly what we 
wanted to prove.
\end{pf*}

\begin{pf*}{Proof of Proposition \r{pr7}}
We first prove the proposition for fixed points using 
Theorem
\r{th3a}.  Let $M$ be the manifold in $J^0(\reals^n, 
\reals^n)$
defined by $f(x) = x$; then $x$ is a fixed point of $f$ if 
and only if
$j^0 f(x)$ lies in $M$.  Notice that $M$ has codimension 
$n$, and
projection onto the first $n$ coordinates is a submersion 
on $M$.  Let
$Z$ be the set of 1-jets in $M \times 
\widehat{P}^1(\reals^n, \reals^n)$
for which $Df$ has an eigenvalue with absolute value $1$.  
Then if $f$
has a nonhyperbolic fixed point, $j^1 f(x)$ must lie in 
$Z$ for some
$x$.  We will be done once we show that $Z$ is a zero set 
in $M \times
\widehat{P}^1(\reals^n, \reals^n)$.  By Lemma \r{le2}, we 
need only show
that the set $S$ of $n \times n$ matrices with an 
eigenvalue on the
unit circle has measure zero (in the space of $n \times n$ 
matrices,
which is isomorphic to $\widehat{P}^1(\reals^n, 
\reals^n)$).  Observe that
every ray from the origin meets $S$ in at most $n$ points, 
because
multiplying every entry of a matrix by a constant factor 
multiplies
its eigenvalues by the same factor.  Hence $Z$ is a zero 
set as
claimed.

For orbits of period $p > 1$, the condition for 
nonhyperbolicity
depends on the 1-jet of $f$ at $p$ different points.  Thus 
to apply
here, Theorem \r{th3a} would have to be generalized to 
$p$-tuples of
jets.  Rather than do this in general, we will prove 
Proposition
\r{pr7} directly from Lemma \r{le3}.  Let us determine 
what conditions
are necessary for a set of polynomial functions $g_1, g_2, 
\ldots, g_q\colon\
\reals^n \to \reals^n$ to span a probe.  For a given $C^1$ 
function
$f\colon\ \reals^n \to \reals^n$ and $\lm \in \reals^q$, let
\[
f_\lm = f + \sum_{i = 1}^q \lm_i g_i.
\]
We must show that for almost every $\lm$, all periodic 
points of
$f_\lm$ with period $p$ are hyperbolic.

Let $x_1, x_2, \ldots, x_p$ denote elements of $\reals^n$, 
and let $A
\sbs \reals^{np}$ be the set of all $(x_1, x_2, \ldots, 
x_p)$ for
which $x_i \neq x_j$ when $i \neq j$.  Consider the 
function $F =
(G,H)$, where $G{:}\  \reals^q \times A \to \reals^{np}$ 
and $H{:}\  \reals^q
\times A \to \reals^{n^2 p}$ are defined by
\[
G(\lm;x_1,x_2,\ldots,x_p) = (f_\lm(x_1) - x_2, f_\lm(x_2) 
- x_3,
\ldots, f_\lm(x_p) - x_1),
\]
\[
H(\lm;x_1,x_2,\ldots,x_p) = (Df_\lm(x_1), Df_\lm(x_2),
\ldots, Df_\lm(x_p)).
\]
(Essentially $F$ consists of the 1-jets of $f$ at $x_1, 
\ldots, x_p$,
except that $G$ projects the 0-jets onto a subspace.)

For a given $\lm$, if $x_1$ is a point of period $p$ for 
$f_\lm$, then
there is a corresponding point $(x_1,\ldots,x_p) \in A$ at 
which $G =
0$.  We then let $M = \{0\}$ in applying Lemma \r{le3}.  
If in
addition $x_1$ is nonhyperbolic, then the matrix 
$\prod_{i=1}^p
Df_\lm(x_i)$ has an eigenvalue on the unit circle.  That is,
$H(\lm;x_1,\ldots,x_p)$ lies in the set $S$ given by
\[
S = \left\{(M_1, \ldots, M_p) \in \reals^{n^2p} : 
\prod_{i=1}^p M_i \text{
has an eigenvalue on the unit circle}\right\},
\]
where $M_1, \ldots, M_p$ denote $n \times n$ matrices.  As 
in our
previous argument for fixed points, $S$ has measure zero 
because every
ray from the origin in $\reals^{n^2 p}$ intersects $S$ in 
at most $n$
points.  Thus we let $Z = \{0\} \times S$ in applying 
Lemma \r{le3}.

We will be done if we can show that $G$ and $H$ satisfy the
hypotheses of Lemma \r{le3}.  Now $G$ and $H$ satisfy the
differentiability hypothesis of Lemma \r{le3} because 
$f_\lm$ is $C^1$
as a function of $x$ and $C^\infty$ as a function of 
$\lm$.  To verify
the transversality hypothesis, we will show that for all 
$(\lm; x_1,
\ldots, x_p) \in \reals^q \times A$, the derivative of $F$ 
with
respect to $\lm$ has full rank. Since $F$ is a linear 
function of
$\lm$, we simply want to show, for every $(x_1, \ldots, 
x_p) \in A$,
that $F$ is onto as a function of $\lm$.  Recall that 
$f_\lm = f +
\sum \lm_i g_i$, and observe that whether or not $F$ is 
onto is
independent of $f$.  We have thus reduced the problem to 
one of 
polynomial interpolation; we need only show there exists a 
finite-dimensional vector space $P$ of polynomial 
functions from $\reals^n
\to \reals^n$ such that for any $p$ distinct points $x_1, 
\ldots, x_p
\in \reals^n$ and any prescribed values for the 1-jet of a 
function
at the $p$ points, there exists a function in $P$ whose 
1-jet takes on
the prescribed values at the prescribed points.

We claim that the polynomials of degree at most $2p - 1$ 
have the
above interpolation property.  We are referring to 
polynomial
functions from $\reals^n \to \reals^n$, but the 
interpolation can be
done separately for each coordinate in the range, so for 
simplicity we
consider polynomials from $\reals^n \to \reals$.  Given 
distinct
points $x_1, \ldots, x_p \in \reals^n$, consider the 
polynomials
\[
P_j(x) = \prod_{\scriptstyle i =1 \atop \scriptstyle i
\scriptstyle \neq\scriptstyle j}^p |x -x_i|^2\]
for $j = 1, \ldots, p$.  Each $P_j$ has degree $2p - 2$ 
and is zero
at every $x_i$ except for $x_j$, where it is nonzero.  
Thus a
suitable linear combination of the $P_j$ can take on any 
prescribed
values at $x_1, \ldots, x_p$.  Next let $P_{jk}(x)$ be the 
$k$th
coordinate of the function $x \mapsto P_j(x) (x - x_j)$ 
for $k = 1,
\ldots, n$.  Each $P_{jk}$ has degree $2p - 1$, and both 
$P_{jk}$ and
its first partial derivatives are all zero at every $x_i$, 
except
that the $k$th partial derivative of $P_{jk}$ is nonzero 
at $x_j$.
Then given a linear combination of the $P_j$ which takes on
prescribed values at $x_1, \ldots, x_p$, adding a linear 
combination
of the $P_{jk}$ will not change these values, and a 
suitable linear
combination of the $P_{jk}$ can be added to change the 
first partial
derivatives at $x_1, \ldots, x_p$ to any prescribed 
values.  This
completes the proof.
\end{pf*}

\begin{pf*}{Proof of Proposition \r{pr8}}
There are two main tasks involved in this proof.  First, 
we must
formulate conditions on the 3-jet of $f$ which must be 
satisfied if
$f$ has an atypical (in the sense of violating one of the 
conditions
in Proposition \r{pr8}) Andronov-H\"{o}pf bifurcation.  
Second, we
must show that the set $Z$ of 3-jets which satisfy these 
conditions
satisfies the hypotheses of Theorem \r{th3a}.

The manifold $M$ which will contain the 2-jet of every 
3-jet in $Z$
consists of those 2-jets which satisfy the following two 
conditions:
\begin{enumerate}
\item[{\ns (a)}] $f = 0$.
\item[{\ns (b)}] $D_x f$ has zero trace and positive 
determinant.
\end{enumerate}
(Of course, these conditions really depend only on the 
1-jet.)
Condition (b) is equivalent to the condition that $D_x f$ 
has nonzero,
pure imaginary eigenvalues.  Notice that condition (a) 
defines a codimension 2
manifold, and adding condition (b) makes $M$ have 
codimension 3.

Condition (i) of Proposition \r{pr8} follows immediately 
from the
implicit function theorem, since the determinant of $D_x 
f$ is nonzero
at the bifurcation point $(\mu_0, x_0)$.  For condition 
(ii) of
Proposition \r{pr8} to hold it suffices that the 
eigenvalues of $D_x
f$ at the fixed point $(\mu, x(\mu))$ have negative real 
parts for
$\mu$ on one side of $\mu_0$ and positive reals parts for 
$\mu$ on the
other side.  Since the eigenvalues of $D_x f$ are complex 
conjugates
in a neighborhood of $(\mu_0, x_0)$, each one has real 
part equal to
half the trace of $D_x f$.  Thus if condition (ii) of 
Proposition
\r{pr8} fails, the following condition must hold:
\begin{enumerate}
\item[{\ns (c)}] The trace of $D_x f(\mu, x(\mu))$ has
  $\mu$-derivative zero at $\mu_0$.
\end{enumerate}
The derivative of this trace depends on the 2-jet of $f$ 
at $(\mu_0,
x_0)$ and on $x\pr(\mu_0)$, which in turn depends on the 
1-jet of $f$.

We wish to show that the set of 2-jets in $M$ for which 
condition (c)
holds is a zero set in $M$.  Since $M$ depends only on the 
1-jet of
$f$, by Lemma \r{le2} it suffices to fix the 1-jet and 
show that as
the second derivatives in the 2-jet vary, condition (c) 
fails almost
everywhere.  Notice that fixing the 1-jet also fixes 
$x\pr(\mu_0)$.
Let the coordinates of $f$ be $(g,h)$ and the coordinates 
of $x$ be
$(y,z)$.  Then condition (c) can be written as
\[
g_{\mu y} + y\pr(\mu_0) g_{yy} + z\pr(\mu_0) g_{yz} + 
h_{\mu z} +
y\pr(\mu_0) h_{yz} + z\pr(\mu_0) h_{zz} = 0,
\]
where the partial derivatives are evaluated at $(\mu_0, 
x_0)$.  As the
second derivatives vary over all real numbers, the above 
equation
holds only on a set of measure zero (a codimension 1 
subspace, in
fact).

We have shown that condition (c) holds only on a zero set 
in $M$.  By
another application of Lemma \r{le2}, the set of 3-jets 
which satisfy
condition (c) is a zero set in $M \times 
\widehat{P}^3(\reals^2,
\reals^2)$.

If conditions (a) and (b) hold while (c) fails, then there 
is a
condition (d) that the 3-jet of $f$ must satisfy in order 
for
condition (iii) of Proposition \r{pr8} to fail.  If 
coordinates
$(u,v)$ are chosen in such a way that
\[
D_x f(\mu_0, x_0) = \lt[ \ba{cc} 0 & -\om \\ \om & 0 \ea 
\rt]
\]
(where $\om$ is the square root of the determinant of $D_x 
f)$, and
$g$ and $h$ are the components of $f$ in this coordinate 
system (this
is different from the definition of $g$ and $h$ above), 
then condition
(d) can be written as
\bea
&\om(g_{uuu} + g_{uvv} + h_{uuv} + h_{vvv}) + 
g_{uv}(g_{uu} + g_{vv})\\
&\qquad  -
h_{uv}(h_{uu} + h_{vv}) - g_{uu} h_{uu} + g_{vv} h_{vv} = 0.
\eea
See \cite{MM} for a detailed derivation of this condition, 
or
\cite{GH} for a more expository discussion of this problem.

Notice that given condition (b), $D_x f$ can be put into 
antisymmetric
form by a linear change of coordinates depending only on 
the 1-jet of
$f$ at $(\mu_0, x_0)$, and further $\om$ is nonzero and 
depends only
on the 1-jet of $f$.  Writing condition (d) in terms of 
the original
coordinates would be tedious; instead we employ Lemma 
\r{le2} again by
fixing the 2-jet of $f$ and letting its third derivatives 
vary.  With
the 2-jet fixed, the above change of coordinates is fixed, 
and induces
a change of coordinates on the space 
$\widehat{P}^3(\reals^2, \reals^2)$.
In terms of the new coordinates, condition (d) determines a
codimension 1 hyperplane in $\widehat{P}^3(\reals^2, 
\reals^2)$, and in
particular the set on which it is satisfied has measure 
zero.
Therefore by Lemma \r{le2}, the set of all 3-jets which 
satisfy
condition (d) is a zero set in $M \times 
\widehat{P}^3(\reals^2,
\reals^2)$.

To summarize, we have shown that in order for the 
conditions given in
Proposition \r{pr8} to fail for a given one-parameter 
family of vector
fields $f$, there must be a point $(\mu_0, x_0)$ at which 
conditions
(a), (b), and at least one of (c) and (d) hold.  We have 
shown that
the manifold $M \sbs J^2(\reals^2, \reals^2)$ defined by 
conditions
(a) and (b) satisfies the hypotheses of Theorem \r{th3a}, 
and that the
set $Z \sbs M \times \widehat{P}^3(\reals^2, \reals^2)$ on 
which at least
one of conditions (c) and (d) holds is a zero set in this 
manifold.
Therefore by Theorem \r{th3a}, for almost every $f$ in $C^k$
the conditions given in Proposition \r{pr8} hold.
\end{pf*}

\begin{pf*}{Proof of Lemma \r{le3}}
We assume without loss of generality that $Z$ is a Borel 
set; then so
is $F^{-1}(Z)$.  Let $\pi_1{:}\ \reals^q \times \reals^n 
\to \reals^q$
be projection onto the first $q$ coordinates.  We wish to 
show that
$\pi_1(F^{-1}(Z))$ has measure zero.  Let $L = G^{-1}(M)$, 
and for
$x \in A$ let $L_x = L \cap (\reals^q \times \{x\})$ be the
``\<$x$-slice'' of $L$. By the transversality hypothesis, 
$G(\cdot,x)$
is transversal to $M$, and thus $L \sbs \reals^q \times 
\reals^n$ and
$L_x \sbs \reals^q$ are manifolds with the same 
codimension, $n$, as
$M \sbs \reals^m$ (see p. 28 of \cite{GP}).

Since $L$ has dimension $q$, away from its critical points 
$\pi_1|_L$
is locally a diffeomorphism.  We will show that 
$F^{-1}(Z)$ is a zero
set in $L$; then since zero sets map to zero sets under
diffeomorphisms, $\pi_1(F^{-1}(Z))$ consists of a zero set 
plus
possibly some critical values of $\pi_1|_L$.  By the Sard 
theorem
\cite{Sar}, the critical values of $\pi_1|_L$ have measure 
zero, and
hence $\pi_1(F^{-1}(Z))$ has measure zero as desired.

To show that $F^{-1}(Z)$ is a zero set in $L$, we first 
show for all
$x \in A$ that $F^{-1}(Z) \cap L_x$ is a zero set in 
$L_x$.  Since
$F(\cdot, x)$ is transversal to $M \times \{y\}$ for all 
$y \in
\reals^s$, and the tangent space $T_\lm L_x$ is the 
inverse image
of $T_{G(\lm,x)} M$ under $D_\lm G$, and both tangent 
spaces have the
same codimension, it follows that $D_\lm F$ maps $T_\lm 
L_x$ onto
$T_{G(\lm,x)} M \times \reals^s$ for all $(\lm,x) \in B 
\times A$.
In other words, $F(\cdot,x)$ is a submersion from $L_x$ to 
$M \times
\reals^s$ for all $x \in A$.  Since $Z$ is a zero set in 
$M \times
\reals^s$, its preimage $F^{-1}(Z) \cap L_x$ is a zero set 
in $L_x$
as claimed.

It remains only to show that the partition $\{L_x\}$ of 
$L$ satisfies
the hypotheses of Lemma \r{le2}.  Let $\pi_2 {:}\  
\reals^q \times
\reals^n \to \reals^n$ be projection onto the last $n$ 
coordinates.
For each $(\lm,x) \in L$, the kernel of $\pi_2$ in 
$T_{(\lm,x)} L$ is
just $T_\lm L_x$.  Since the former space has dimension 
$q$ and the
latter space has dimension $q - n$, it follows that 
$\pi_2$ has rank
$n$ on $T_{(\lm,x)} L$.  Thus $\pi_2|_L$ is a submersion, 
which
implies (see p. 20 of \cite{GP}) that near every point in 
$L$ there is
a local $C^1$ coordinate system on $L$ whose last $n$ 
coordinates are
the same as those of $x$.  The slices $L_x$ of $L$ are 
parallel
hyperplanes in such a coordinate system, and therefore the 
hypotheses
of Lemma \r{le2} are satisfied.
\end{pf*}

\section{Extensions of prevalence} \label{se5}

In this article we have proposed sufficient conditions for 
a property
to be said to be true ``almost everywhere'', in a 
measure-theoretic
sense, on complete metric linear spaces.  In other 
contexts more
general definitions may be appropriate.  For instance, the 
concepts
of shyness and prevalence can be extended from vector 
spaces to larger
classes of topological groups \cite{Myc}.

We have concentrated thusfar on extending the notions of 
``measure
zero'' and ``almost every'' to infinite-dimensional 
spaces.  We now
briefly consider some ways to characterize sets which are 
neither shy
nor prevalent in an infinite-dimensional vector space $V$.
\begin{definition} \label{df10} \ns
Let $P$ be the set of compactly supported probability 
measures on
the Borel sets of $V$.  The {\em lower density\/} 
$\rho^-(S)$ of a
Borel set $S \sbs V$ is defined to be
\[
\rho^-(S) = \sup_{\mu \in P} \inf_{v \in V} \mu(S + v).
\]
The {\em upper density\/} $\rho^+(S)$ is given by
\[
\rho^+(S) = \inf_{\mu \in P} \sup_{v \in V} \mu(S + v).
\]
If $\rho^-(S) = \rho^+(S)$, then we call this number the 
{\em
relative prevalence\/} of $S$.
\end{definition}

One can show that for all $\mu, \nu \in P$,
\[
\inf_{v \in V} \mu(S + v) \leq \inf_{v \in V} \mu*\nu(S + v)
\leq \sup_{v \in V} \mu*\nu(S + v) \leq \sup_{v \in V} 
\nu(S + v),
\]
and thus $0 \leq \rho^-(S) \leq \rho^+(S) \leq 1$ for all 
Borel sets
$S$.  It follows that a shy set has relative prevalence 
zero and a
prevalent set has relative prevalence one.  However, sets 
with
relative prevalence zero need not be shy; all bounded sets 
have
relative prevalence zero, for example.

In $\reals^n$, having positive lower density is a much 
stronger
condition on a set than having positive Lebesgue measure.  
The
following weaker conditions give a closer analogue to 
positive
measure.
\begin{definition} \label{df11} \ns
A measure $\mu$ is said to {\em observe\/} a Borel set $S 
\sbs V$ if
$\mu$ is finite and $\mu(S + v) > 0$ for all $v \in V$.  A 
Borel set
$S \sbs V$ is called {\em observable\/} if there is a 
measure which
observes $S$, and is called {\em substantial\/} if it is 
observed by
a compactly supported measure.  More generally, a subset 
of $V$ is
observable (resp. substantial) if it contains an 
observable (resp.
substantial) Borel set.
\end{definition}

Every set with positive lower density is then substantial, 
and every
substantial set is observable.  As in Fact \r{fa3}, if 
$\mu$ observes
a Borel set $S$ then so does $\mu*\nu$ for any finite 
measure $\nu$.
It follows that an observable set is not shy.  In 
$\reals^n$, it
follows as in Fact \r{fa6} that a set is observable if and 
only if it
contains a set of positive Lebesgue measure.  In a 
separable space
every open set is observable; given a countable dense 
sequence
$\{x_n\}$, the measure consisting of a mass of magnitude 
$2^{-n}$ at
each $x_n$ observes each open set.

\section*{Acknowledgments}

This paper has been greatly enhanced by the input of many 
people, only
a few of whom we can acknowledge individually.  We thank 
P. M.
Fitzpatrick for pointing out the existence of \cite{QS}, a 
paper on a
Sard-type theorem for Fredholm operators, in which the 
concept of
``0-preconull'' is introduced. This concept is similar in 
principle to
a local version of shyness using measures which have 
finite-dimensional 
support.  We thank A. Kagan for bringing to our attention
the results of Sudakov on the nonexistence of 
quasi-invariant
measures, and
J. Milnor for suggesting the two examples from complex
dynamics in the introduction.  We thank Yu.\ Il\cprime 
yashenko for helping to
simplify our proof of Fact \r{fa8} and helping to 
formulate the
approach taken in \S \r{se4}.  We thank C. Pugh for his 
extensive
comments, in particular suggesting the line of reasoning 
used for the
proof of Lemma \r{le3}.  We thank H. E. Nusse and M. Barge 
for their
helpful comments.

The first author is grateful to the ONT Postdoctoral 
Fellowship
Program (administered by the American Society for 
Engineering
Education) for facilitating his work on this paper.

\end{document}